\documentclass[12pt]{article}
\usepackage{amsmath}
\usepackage{amsfonts,euscript,amssymb,amsmath,amsthm}
\usepackage{amstext}
\usepackage[english]{babel}

\newtheorem{theorem}{Theorem}

\newtheorem{lemma}{Lemma}
\newtheorem{corollary}{Corollary}

\newenvironment{pfc}{\par\noindent\textbf{Proof of the statement c).}}{\hfill$\square$}
\newenvironment{pfb}{\par\noindent\textbf{Proof of the statement b).}}{\hfill$\square$}
\newenvironment{pf0}{\par\noindent}{\hfill$\square$}

\newcommand{\bbQ}{\mathbb{Q}}

\newcommand{\bbN}{\mathbb{N}}
\newcommand{\bbZ}{\mathbb{Z}}
\newcommand{\Image}{\mathop{\rm Im}}

\newcommand{\Real}{\mathop{\rm Re}}

\begin{document}

\begin{center}

{\bf \large On the critical line zeros of $L$~---functions attached
to automorphic cusp forms.}

I.S. Rezvyakova \footnote{This work was supported by grant RFBR
11-01-00759a and
 12-01-31165a. }

\end{center}

\title{On the critical line zeros of $L$~---functions attached to
automorphic cusp forms.}







\begin{center}
{\bf \S 1. Introduction. Statement of the main result}
\end{center}


One of the most interesting questions in the theory of the Riemann
zeta-function is the Riemann hypothesis which asserts that all
non-trivial zeros of the Riemann zeta-function lie on the critical
line. The Riemann hypothesis is not yet proved nor disproved. In
1942, Atle Selberg \cite{Selberg1942} showed that a positive
proportion of non-trivial zeros of the Riemann zeta-function lie on
the critical line (for the numerical estimates
see works~\cite{Levinson}--\cite{Conrey3}).
The same type result holds also for Dirichlet $L$~- functions
(\cite{Zhuravlev2}). In 1989, A.~Selberg in his report at the
conference in Amalfi conjectured that all functions from Selberg
class $S$ (which have decomposition in the Euler product and
functional equation of the Riemann type as the necessary conditions)
satisfy an analogue of the Riemann  hypothesis
(see~\cite{Selberg_Amalfi}).

The Riemann zeta-function and Dirichlet $L$~- functions are
functions of degree one (for the definition of the degree, which is
a characteristic of a functional equation,
see~\cite{Selberg_Amalfi}). In 1983, J.L.~Hafner proved an analogue
of Selberg's theorem for a function of degree two.
In his papers~\cite{Hafner1}, \cite{Hafner2} for $L$~- series, whose
coefficients are attached to those of holomorphic cusp forms
(modular forms) or non-holomorphic cusp forms (Maass wave forms) for
the full modular group with trivial character, the result on
positivity of proportion of non-trivial zeros lying on the critical
line is established (see the necessary definitions further in the
text and in~\cite{Iwaniec_book}). This work is a continuation of
~\cite{Hafner1}, \cite{Hafner2}.
Here we obtain an analogue of Selberg's theorem for $L$~- functions
attached to automorphic cusp forms with respect to the Hecke
congruence subgroup $\Gamma_0(D)$ with arbitrary integral weight
$k\ge 1$ (theorem~\ref{th}). First, we recall some definitions and
notation.

Suppose that $f(z)$ is an automorphic cusp form of integral weight
$k\ge 1$ for the group $\Gamma_0 (D)$ with a character $\chi$ modulo
$D$ (briefly we write this as $f\in S_k(\Gamma_0 (D), \chi)$), where
$\Gamma_0(D)$ is a subgroup of $\mbox{SL}_2 (\bbZ)$ that consists of
all the matrices
$\gamma = \bigl( \begin{smallmatrix} a & b \\
c & d
\end{smallmatrix}\bigr)$ satisfying the condition $ c \equiv
0 (\!\!\!\mod D)$). In other words, let $f$ be a holomorphic
function on the upper half-plane, such that for every element
$\gamma\in \Gamma_0(D)$ the following relation is fulfilled:
\begin{equation}\label{eq1}
f(\gamma z) = \chi (\gamma) (cz +d)^{k} f(z), \quad \text{ where }
\quad \chi (\gamma) = \chi (d), \quad \gamma z = \frac{az+b}{cz+d},
\end{equation}
and also that $f$ vanishes at every cusp of the group $\Gamma_0
(D)$. This entails that
$f(z)$ has an expansion
$$
f(z)=\sum\limits_{n=1}^{+\infty} a(n) e^{2\pi i nz} \quad \text{for}
\quad \Real z >0.
$$
Next, assume that $f$ is not identically zero and is an
eigenfunction of all the Hecke operators $T_n$ for $n=1,2,\ldots$,
where
\begin{equation}\label{eq0}
T_n f (z) = \frac{1}{n} \sum\limits_{ad = n} \chi (a) a^k
\sum\limits_{0 \le b <d} f \left( \frac{az+b}{d}\right).
\end{equation}
Without loss of generality we may assume that $a(1)=1$. From the
properties of the Hecke operators, the equality $T_{n} f = a(n) f$
follows for every positive integer $n$, and also if $\Real s
>1$, then for Dirichlet series
\begin{equation}\label{eq3}
L(s)= L_{f}(s)=\sum\limits_{n=1}^{+\infty} \frac{r(n)}{n^{s}}
\end{equation}
with
\begin{equation}\label{eq2}
r(n) = a(n) n^{\frac{1-k}{2}},
\end{equation}
the identity
$$
L(s)= \prod_{p} \left( 1- \frac{r(p)}{p^s} + \frac{\chi
(p)}{p^{2s}}\right)^{-1}
$$
holds (here the product is carried over
all consecutive prime numbers). For the normalized coefficients
$r(n)$ of the automorphic (for a congruence subgroup) cusp form of
an integral weight, which is also an eigenfunction of all the Hecke
operators, the estimate
$$
|r(n)| \le \tau(n)
$$
is valid, where $\tau(n)$ is the number of divisors of $n$. This
inequality was previously known as the Ramanujan~- Petersson
conjecture until its truth was proved in~\cite{Deligne1974},
\cite{Del_Serre1974}. The function $L(s)$ satisfies the following
functional equation (see~\cite{Iwaniec_book}, \S 6.7, \S 7.2 ):
$$
\Lambda(s) = \theta \cdot \overline{\Lambda(1-\overline{s})},
$$
where $\Lambda(s)$ is an entire function,
$$
\Lambda(s) = \left(\frac{2\pi}{\sqrt{D}}\right)^{-s-\frac{k-1}{2}}
\Gamma\left( s+\frac{k-1}{2}\right) L(s),
$$
and $|\theta|=1$, or, more precisely , $\theta = i^{k}
\overline{\eta}$, where $\eta$ is the eigenvalue of the operator
$\overline{W}$,
$$
\overline{W} f = (-1)^{k} D^{-k/2} z^{-k} \overline{f\left(
1/D\overline{z}\right)}, \quad \overline{W} f = \eta f.
$$
The function $L(s)$ has zeros at $s=-\frac{k-1}{2},
-\frac{k-1}{2}-1, -\frac{k-1}{2}-2, \ldots$, which correspond to the
poles of $\Gamma(s+\frac{k-1}{2})$ and are called ``trivial zeros''.
The remaining zeros lie in the strip $0 \le \Real s \le 1$ and are
called ``non-trivial''. $L(s)$ is a function of degree two and
belongs to Selberg class. Therefore, an analogue of the Riemann
hypothesis exists for this function which claims that all its
non-trivial zeros lie on the critical line $\Real s = \frac12$. In
this work we prove the following theorem.
\begin{theorem}\label{th}
Suppose that $f(z) = \sum\limits_{n=1}^{+\infty} a(n) e^{2\pi i nz}$
is an automorphic cusp form of integral weight $k \ge 1$ for the
group $\Gamma_0(D)$ with character $\chi$ modulo $D$, which is also
an eigenfunction of all the Hecke operators $T_n$ for $n=1, 2,
\ldots$. Define $L(s) = L_{f} (s)$
by the equalities (\ref{eq3}) and (\ref{eq2}), and set $N_0(T)$ to
be
the number of odd order zeros of $L(s)$ on the interval $\{s=\frac12
+ it, \, 0 < t \le T\}$. Then $N_0(T) \ge c T\ln T$ with some
constant $c>0$.
\end{theorem}

Notice that if $N(T)$ is the number of zeros of $L(s)$ in the
rectangle $\{ s \mid 0\le \Real s \le 1, \; 0< \Image s \le T \}$,
then the asymptotic formula holds
$$
N(T) = \frac{T}{\pi} \ln \frac{T}{\pi c_1} + O(\ln T) \quad
\text{when} \quad T \to +\infty;
$$
thus we have
\begin{corollary}
A positive proportion of non-trivial zeros of $L_{f}(s)$ lie on the
critical line.
\end{corollary}
An example of $L(s)$ considered in this work is a Hecke $L$~-
function with complex class group character on ideals of imaginary
quadratic field $\bbQ (\sqrt{-D})$ (since there is a corresponding
form $f$ from  $S_1(\Gamma_0 (D), \chi_D)$).

\begin{center}
{\bf \S 2. Main and auxiliary statements}
\end{center}

Hereafter we suppose that $k$ (the weight of the form $f$) and $D$
(the level of the form $f$) are fixed numbers.

The main idea of the proof belongs to A.~Selberg which is served by
introducing a ``mollifier''. Define numbers $\alpha(\nu)$ by the
equality
\begin{equation}\label{eq4}
\sum\limits_{\nu =1}^{+\infty} \alpha(\nu) \nu^{-s} = \prod_{p >
256} \left( 1-\frac{r(p)}{2 p^{s}} \right).
\end{equation}
Suppose that $X\ge 3$ and set
\begin{equation}\label{eq5}
\beta(\nu ) = \alpha(\nu ) \left( 1-\frac{\ln \nu}{\ln
X}\right)^{+}, \quad \text{where} \quad x^{+}= \max (0, x).
\end{equation}
Let us define a mollifier $\varphi(s)$ by the formula
$$
\varphi(s) = \sum\limits_{\nu =1}^{+\infty} \beta(\nu) \nu^{-s}.
$$

Suppose $\delta >0$,
\begin{equation}\label{eq6}
\mathfrak{F}(t)= \frac{1}{\sqrt{2\pi}} \Lambda\left( \frac12 + it
\right) \left| \varphi\left( \frac12 +it\right)\right|^{2}
\exp\left( \left(\frac{\pi}{2} -\delta\right) t\right).
\end{equation}
The functional equation for $L(s)$ yields
that $\theta^{-1/2} \mathfrak{F}(t)$ is real-valued
for real $t$. Observe also,
that odd order zeros of $\theta^{-1/2} \mathfrak{F}(t)$ are those of
$L (1/2 +it)$.

For $0<h_1<1$, put
$$
I_1(t) = \int\limits_{-h_1}^{h_1} |\theta^{-1/2} \mathfrak{F}(t+u)|
du = \int\limits_{-h_1}^{h_1} |\mathfrak{F}(t+u)| du,
$$
$$
I_2(t) = \left|\int\limits_{-h_1}^{h_1} \theta^{-1/2}
\mathfrak{F}(t+u) du \right| = \left|\int\limits_{-h_1}^{h_1}
\mathfrak{F}(t+u) du \right|.
$$
Let $T>1$ be a sufficiently large number. Define $E_1$ as a set of
points $t\in (1, T)$, such that for $t\in E_1$ the inequality
$$
I_1(t) > I_2(t)
$$
holds. By $E_2$ denote the complementary set to $E_1$, i.e., the set
of all points $t\in (1, T)$ with
$$
I_1(t) = I_2(t).
$$
If $\mu(E_1)$ denotes the measure of the set $E_1$, then the number
of odd order zeros of the function $\theta^{-1/2}\mathfrak{F}(t)$ on
the interval $(0, T)$ is not less than $\frac{\mu(E_1)}{2 h_1} -1$
(see~\cite{Selberg1942} or \cite{Voronin_Karatsuba}).

Now the main statement of this work is a consequence of the
following theorem.

\begin{theorem}\label{th1}
Suppose that $0<\delta < \delta_0$, where $\delta_0 <1/10$ is some
small positive number, $0 < h_1 <1$, $X \ge 3$ and $\delta X^{86}
e^{1/h_1} \le 1$. Then the following estimates are valid:

a) $\int\limits_{1}^{\delta^{-1}} I_1 (t) dt \gg \delta^{-1} h_1$,

b) $\int\limits_{-\infty}^{+\infty} I_1^{2} (t) dt \ll \delta^{-1}
\dfrac{h_1^{2} \ln \delta^{-1}}{\ln X}$,

c) $\int\limits_{-\infty}^{+\infty} I_2^{2} (t) dt \ll \delta^{-1}
\dfrac{h_1}{\ln X}$, where the constants implied in Vinogradov's
signs $\ll, \gg$ are absolute.
\end{theorem}

To derive
the main result, set
in theorem~\ref{th1}
$$\delta
= T^{-1}, \quad X = T^{1/100}, \quad h_1 = \frac{A}{\ln X},
$$ where $A$ is a sufficiently large positive constant,
and write the following chain
of relations
$$
I_3 = \int\limits_{1}^{T} I_1(t) dt = \int\limits_{E_1} I_1(t) dt +
\int\limits_{E_2} I_1(t) dt = \int\limits_{E_1} I_1(t) dt +
\int\limits_{E_2} I_2(t) dt \le I_1 + I_2,
$$
where
\begin{equation*}
\begin{split}
I_1 &= \int\limits_{E_1} I_1(t) dt \le \left({\mu(E_1)}
\right)^{1/2} \left(
\int\limits_{-\infty}^{+\infty} I_1^{2} (t) dt \right)^{1/2}, \\
I_2 &= \int\limits_{1}^{T} I_2 (t) dt \le T^{1/2} \left(
\int\limits_{-\infty}^{+\infty} I_2^{2} (t) dt\right)^{1/2}.
\end{split}
\end{equation*}
The result c) of theorem~\ref{th1} entails
$$
I_2 \ll T \left( \dfrac{h_1}{\ln X}\right)^{1/2}.
$$
Hence, by virtue of the
estimate a), for sufficiently large
$A$ we have
$$
I_3 \ge 2 I_2.
$$
Therefore,  $I_1 \ge I_3/2$, and thus the relations a) and b) of
theorem~\ref{th1} imply the
estimates
$$
T h_1 \ll I_3 \le 2 I_1 \ll (\mu E_1)^{1/2} \left( T h_1^{2}
\dfrac{\ln T}{\ln X}\right)^{1/2},
$$
or
$$
\mu(E_1) \gg T.
$$
Thereby, the number of odd order zeros of
$\theta^{-1/2}\mathfrak{F}(t)$ (that is of $L(1/2+it)$) on the
interval $0 < t \le T$ is estimated from below by the quantity of
order
$$
\mu(E_1) h_1^{-1} \gg T \ln X \gg T \ln T.
$$
Now it only remains
to establish the statement of theorem~\ref{th1}.

Proof of the assertion a) of theorem~\ref{th1} repeats the proof of
a similar relation while considering the Riemann zeta-function
instead of $L(s)$ (see~\cite{Voronin_Karatsuba}, \S 6.3).

To derive
the assertions b) and c) of theorem~\ref{th1} we employ
the following auxiliary lemma.
\begin{lemma}\label{l1}
Suppose that $\mathfrak{F}(y)$ is given by the formula (\ref{eq6}),
and
\begin{equation}\label{eq7}
\begin{split}
G(y) = \left| \sum\limits_{n,\nu_1, \nu_2 \in \bbN}
\frac{r(n)\beta(\nu_1)\overline{\beta(\nu_2)}}{\nu_2} \exp \left(
-\frac{2\pi n \nu_1}{\sqrt{D}\nu_2} y (\sin\delta +i\cos
\delta)\right) \right|^2.
\end{split}
\end{equation}
Then
$$
\int\limits_{-\infty}^{+\infty} \left( \int\limits_{-h_1}^{h_1}
|\mathfrak{F}(t+u)| du \right)^2 dt \le 8 h_1^2
\int\limits_{1}^{+\infty}G(y) dy,
$$
and
\begin{equation*}
\begin{split}
\int\limits_{-\infty}^{+\infty}\left|\int\limits_{-h_1}^{h_1}
\mathfrak{F}(t+u) du \right|^2 dt  & \le 8 h_1^2 \int\limits_{1}^{H}
G(y) dy   + 8 \int\limits_{H}^{+\infty} \frac{G(y)}{\ln^2 y} dy,
\end{split}
\end{equation*}
where $H = e^{1/h_1}$.
\end{lemma}

The proof of this result is contained in~\cite{Rezvyakova} (lemma 3
and 4).

Let us formulate
now the main lemmas, from which the statement of theorem~\ref{th1}
will follow easily.

\begin{lemma}\label{l2}
Assume that $G(y)$ is defined by the formula (\ref{eq7}). Then under
the conditions of theorem~\ref{th1}, for $1 \le x \le e^{1/h_1}$ the
estimate
$$
J(x, \theta) = \int\limits_{x}^{+\infty} G(u) u^{-\theta} du \ll
\frac{\delta^{-1}}{\theta x^{\theta} \ln X}
$$
is valid
uniform in $\theta$ from the interval $0 < \theta
\le 1/4$.
\end{lemma}

The core of the proof of
lemma~\ref{l2} relies on
the following two lemmas. Lemma~\ref{l3} is related to an estimation
of a ``diagonal'' term (estimation of Selberg sums), and
lemma~\ref{l4} accordingly to a ``non-diagonal'' term. Let us adopt
further the following notation:
\begin{itemize}
\item $\overline{K}(s) = \overline{K(\overline{s})}$,
\item $p^{\alpha} || m$ means that $p^{\alpha}$ divides $m$, but
$p^{\alpha+1} \nmid m$.
\end{itemize}

\begin{lemma}\label{l3}
Let $0 \le \theta \le 1/4$, and define the sum $S(\theta)$ by the
equality
$$
S(\theta)=  \sum\limits_{\nu_1,\ldots, \nu_4 \le X}
\frac{\beta(\nu_1) \overline{\beta(\nu_2)} \overline{\beta(\nu_3)}
\beta(\nu_4)}{ \nu_2 \nu_4} \left( \frac{q}{\nu_1 \nu_3}
\right)^{1-\theta} K\left( \frac{\nu_1 \nu_4}{q}, 1-\theta\right)
\overline{K}\left( \frac{\nu_2 \nu_3}{q}, 1-\theta\right),
$$
where
$$
q=(\nu_1 \nu_4, \nu_2\nu_3),
$$
\begin{equation*}
\begin{split}
K(m,s) &=  \prod\limits_{p| m}
\left( 1+ \frac{|r^{2}(p)|}{p^{s}} + \frac{|r^{2}(p^2)|}{p^{2s}}+\ldots \right)^{-1} \times \\
&\times \prod_{p^{\alpha} || m} \left( \overline{r(p^{\alpha})} +
\frac{\overline{r(p^{\alpha+1})}r(p)}{p^{s}} +
\frac{\overline{r(p^{\alpha+2})}r(p^2)}{p^{2s}}+\ldots \right).
\end{split}
\end{equation*}
Then the estimate
$$
S(\theta) \ll \frac{X^{2\theta}}{\ln X}
$$
holds uniformly in $\theta$.
\end{lemma}

\begin{lemma}\label{l4}
Let $N\gg 1$, $(m_1, m_2) =1$, $m_1^{8} m_2^{9} \le N$, \; $l \le
N^{10/11}$,
\begin{equation*}
\begin{split}
S = \sum\limits_{n=1}^{N-1} r(n) \overline{r \left( \frac{m_1
n+l}{m_2} \right)}
\end{split}
\end{equation*}
(assuming that
the function $r(\cdot)$ vanishes for non-integral argument). Then
for arbitrary $\varepsilon
>0$ the following estimate is valid:
\begin{equation*}
\begin{split}
S  \ll_{\varepsilon} N^{10/11+\varepsilon} m_1^{8/11} m_2^{-2/11}.
\end{split}
\end{equation*}
\end{lemma}

Lemma~\ref{l4} implies the following
\begin{corollary}\label{cor1}
The function
$$
D_{m_1, m_2} (s,l) = \sum\limits_{n=1}^{+\infty} \frac{r(n)
\overline{r \left( \frac{m_1 n+l}{m_2} \right)}}{(m_1 n + l/2)^{s}}
$$
has an analytic continuation
in the half-plane $\Real s
> 10/11$, and, moreover, for $0<\varepsilon_0 \le 1/11$, in the
region $\Real s \ge 10/11 +\varepsilon_0$ the estimate
$$
D_{m_1, m_2} (s,l) \ll_{\varepsilon_0}
\frac{|s|}{m_1^{\frac{10}{11}+\varepsilon_0} m_2} \left( (m_1^{8}
m_2^{9})^{\frac{1}{11}-\frac{\varepsilon_0}{2}} + l^{\frac{1}{10}-
\frac{11\varepsilon_0}{20}} \right)
$$
holds true.
\end{corollary}

The deduction of
corollary~\ref{cor1} from lemma~\ref{l4} is contained
in~\cite{Rezvyakova}. The main tool for obtaining the statement of
lemma~\ref{l4} is Jutila's variant of circle method (\cite{Jutila},
\cite{Jutila1992}), which we also used to prove similar result
\cite{Rezvyakova} (lemma~7) for coefficients $r(n)$ of automorphic
cusp forms of weight $k=1$.
A modification
in the proof of the statement for an arbitrary weight $k$ is that
one has to use
the following analogue of lemma~10 in~\cite{Rezvyakova}: {\it let
$q\equiv 0 (\!\!\!\mod D)$, $(a,q)=1$, $a a^{*} \equiv 1( mod q)$,
and let $k(t)$ be a smooth (for example, twice continuously
differentiable) function with compact support on $(0, +\infty)$.
Then
\begin{equation*}
\begin{split}
\sum\limits_{n=1}^{+\infty} r(n) k(n) e^{2\pi i \frac{a}{q}n} =
\overline{\chi (a)} \sum\limits_{n=1}^{+\infty} r(n) e^{-2\pi i
\frac{a^{*}}{q}n} \tilde{k}(n),
\end{split}
\end{equation*}
where $\tilde{k} (n) = \frac{2\pi i^k}{q} \int\limits_{0}^{+\infty}
k(t) J_{k-1} \left( \frac{4\pi \sqrt{nt}}{q}\right) dt$, and
$J_{k-1} (t)$ is Bessel function. }

To establish the above
relation one has to notice for an automorphic cusp form $f \in
S_k(\Gamma_0(D), \chi)$,
$$
f(z) = \sum\limits_{n=1}^{+\infty} a(n) e^{2\pi i nz},
$$
that the equality
\begin{equation}\label{eq8}
f\left( \frac{a}{q} + iw \right) = \frac{\overline{\chi
(a)}}{(-iqw)^{k}} f\left( - \frac{a^{*}}{q} - \frac{1}{iq^2
w}\right), \quad \Real w >0,
\end{equation}
holds, which implies the relation
\begin{equation}\label{eq9}
\begin{split}
& \sum\limits_{n=1}^{+\infty} a(n)  k(n) e^{2\pi i \frac{a}{q} n} =
i^{k-1}\frac{\overline{\chi (a)}}{q^k} \sum\limits_{n=1}^{+\infty}
a(n) e^{-2\pi i \frac{a^{*}}{q} n}\int\limits_{0}^{+\infty} k(t)
\left(\; \int\limits_{c-i\infty}^{c+i\infty} \frac{1}{w^k} e^{-
2\pi\frac{n}{q^2 w}} e^{2\pi tw} dw \right)dt.
\end{split}
\end{equation}
Thus, by virtue of
the identity
$$
\frac{1}{2\pi i} \int\limits_{c-i\infty}^{c+i\infty} \frac{e^{as
-bs^{-1}}}{s^{k}}
 ds = \left(\frac{a}{b}\right)^{\frac{k-1}{2}} J_{k-1}(2\sqrt{ab}),
$$
which is valid for $a,b,c >0$, $k \ge 1$ (see~\cite{Beitman}, \S
7.3.3), we have that the integral enclosed in parenthesis
in the previous relation is equal to
$$
2\pi i q^{k-1} \left(\frac{t}{n}\right)^{\frac{k-1}{2}}J_{k-1}
\left( \frac{4\pi \sqrt{nt}}{q}\right).
$$
To complete the proof of that auxiliary statement one has to set
$k(t) = k_1 (t) t^{\frac{1-k}{2}}$ in the formula (\ref{eq9})
and recall that $r(n) = a(n) n^{\frac{1-k}{2}}$.

\begin{center}
{\bf \S 3. Proof of statements b) and с) of theorem~\ref{th1}}
\end{center}

Let us prove assertions b) and c) of theorem~\ref{th1} using the
lemmas formulated above.
\begin{pfc}
For $H= e^{1/h_1}$, applying lemma~\ref{l2} with $\theta = 1/4$, we
get the estimate
\begin{equation*}
\begin{split}
\int\limits_{1}^{H} G(x) dx & = -\int\limits_{1}^{H} x^{\theta}
\frac{d}{dx} J(x, \theta) dx = \left. -x^{\theta}  J(x, \theta)
\right|_{x=1}^{x=H} + \theta \int\limits_{1}^{H} x^{\theta-1} J(x,
\theta) dx \\
&\ll \frac{\delta^{-1}}{\theta \ln X} +  \frac{\delta^{-1} \ln
H}{\ln X} \ll \frac{\delta^{-1}}{h_1 \ln X}.
\end{split}
\end{equation*}

Similarly, with $\theta = 1/4$ we have
\begin{equation*}
\begin{split}
\int\limits_{H}^{+\infty} \frac{G(x)}{\ln^2 x} dx & =
-\int\limits_{H}^{+\infty} \frac{x^{\theta}}{\ln^2 x}
\frac{d}{dx} J(x, \theta) dx \\
&= \left. -\frac{x^{\theta}}{\ln^2 x}  J(x, \theta)
\right|_{x=H}^{x=+\infty} + \int\limits_{H}^{+\infty} J(x,\theta)
x^{\theta}
\left( \frac{\theta}{x\ln^2 x} -\frac{2}{x\ln^3 x} \right) dx \\
&\ll \frac{\delta^{-1}}{\theta \ln X \ln^2 H} +
\frac{\delta^{-1}}{\ln X \ln H} + \frac{\delta^{-1}}{\theta \ln X
\ln^2 H} \ll \frac{\delta^{-1} h_1 }{\ln X}.
\end{split}
\end{equation*}
Now the statement с) follows from lemma~\ref{l1}.
\end{pfc}

\begin{pfb}
Application of lemma~\ref{l2} with $x=1$, $\theta = \dfrac{1}{\ln
\delta^{-1}}$ gives
$$
\int\limits_{1}^{\delta^{-2}} G(u) du  \ll \int\limits_{1}^{+\infty}
G(u) u^{-1/ \ln \delta^{-1}} du \ll \frac{\delta^{-1} \ln
\delta^{-1}}{\ln X}.
$$
Using formula (\ref{eq7}), we estimate the integral on the interval
$(\delta^{-2},+\infty )$ of $G(u)$ by the following expression
\begin{equation*}
\begin{split}
\int\limits_{\delta^{-2}}^{+\infty} G(u) du & \ll
\sum\limits_{\substack{n_1, n_2 \\ \nu_1,\nu_2, \nu_3, \nu_4}}
\frac{|r(n_1)r(n_2)\beta(\nu_1)\beta(\nu_2)\beta(\nu_3)\beta(\nu_4)|}{\nu_2
\nu_4}\times \\
&\times  \int\limits_{\delta^{-2}}^{+\infty} \exp \left(
-\frac{2\pi}{\sqrt{D}} \left(\frac{n_1 \nu_1}{\nu_2}+ \frac{n_2
\nu_3}{\nu_4} \right) \delta x \right) dx.
\end{split}
\end{equation*}
Let $A=\frac{2\pi}{\sqrt{D}} \left(\frac{n_1 \nu_1}{\nu_2}+
\frac{n_2 \nu_3}{\nu_4} \right) \delta$. From the equality
$$
\int\limits_{\delta^{-2}}^{+\infty}  \exp (-Ax) dx = \frac{e^{-A
\delta^{-2}}}{A}
$$
we obtain the following estimate:
\begin{equation*}
\begin{split}
\int\limits_{\delta^{-2}}^{+\infty} G(u) du  &\ll \delta^{-1} \!
\sum\limits_{\substack{n_1, n_2 \\ \nu_1,\nu_2, \nu_3, \nu_4}}
\frac{|r(n_1)r(n_2)\beta(\nu_1)\beta(\nu_2)\beta(\nu_3)\beta(\nu_4)|}{(n_1
\nu_1 \nu_4 + n_2 \nu_2 \nu_3)} \exp \left( -\frac{2\pi}{\sqrt{D}}
\delta^{-1} \left(\frac{n_1 \nu_1}{\nu_2}+ \frac{n_2 \nu_3}{\nu_4}
\right)
\right) \\
& \ll \delta^{-1} \! \sum\limits_{\substack{n_1, n_2 \\ \nu_1,\nu_2,
\nu_3, \nu_4}}
\frac{|r(n_1)r(n_2)\beta(\nu_1)\beta(\nu_2)\beta(\nu_3)\beta(\nu_4)|}{\sqrt{n_1
n_2 \nu_1 \nu_2 \nu_3 \nu_4}}  \exp \left( -\frac{2\pi}{\sqrt{D}}
\delta^{-1} \left(\frac{n_1 \nu_1}{\nu_2}+ \frac{n_2 \nu_3}{\nu_4}
\right)
\right) \\
&\le \delta^{-1} \left(  \sum\limits_{\nu_1,\nu_2 \le X}
\frac{1}{\sqrt{\nu_1 \nu_2}} \sum\limits_{n} \frac{|r(n)|}{\sqrt{n}}
\exp \left( -\frac{2\pi}{\sqrt{D}} \frac{n \nu_1}{\nu_2} \delta^{-1}
\right) \right)^2.
\end{split}
\end{equation*}
Since $|r(n)| \le \tau(n)$, and $\nu_2 \le X \le \delta^{-1/3}$,
then
$$
\sum\limits_{n=1}^{+\infty} \frac{|r(n)|}{\sqrt{n}} \exp \left(
-\frac{2\pi}{\sqrt{D}} \frac{n \nu_1}{\nu_2} \delta^{-1} \right) \ll
\exp \left( -\sqrt{\delta^{-1}}\right).
$$
Hence,
$$
\int\limits_{\delta^{-2}}^{+\infty} G(u) du  \ll \delta^{-1} X \exp
\left( -\sqrt{\delta^{-1}}\right) \ll 1.
$$
Now the statement b) of theorem~\ref{th1} follows from
lemma~\ref{l1}.
\end{pfb}

\begin{center}
{\bf \S 4. Proof of main lemma~\ref{l2}}
\end{center}
\begin{pf0}

Similarly to the proof in~\cite{Hafner1}, we introduce a
non-negative smooth factor $\phi(u) \in {\bf C}^2 (0, +\infty)$,
having the property
\begin{align*}
\phi (u) = \left\{
\begin{array}{l}
0, \mbox{ если } u\le 1/2 \\
1, \mbox{ если } u \ge 1.
\end{array}
\right.
\end{align*}
Then
\begin{equation}\label{eq10}
J(x, \theta) \le \int\limits_{0}^{+\infty} \phi \left(
\frac{u}{x}\right) G(u) u^{-\theta} du = x^{1-\theta}
\int\limits_{0}^{+\infty} \phi(u) G (ux)u^{-\theta} du.
\end{equation}
We need several auxiliary relations for estimation of the last
integral. For $\Real s>1$, define
$$
\Phi(s,y) =\int\limits_{0}^{+\infty}\phi(u) u^{-s} \exp(-2\pi iyu)
du.
$$
We then have
$$
\Phi^{*} (s) := \int\limits_{0}^{+\infty} \phi' (u) u^{-s+1} du =
(s-1) \Phi (s,0).
$$
Since $\phi'$ is a function with compacts support, then $\Phi^{*}
(s)$ is an entire function. Therefore, $\Phi^{*} (s)$ is the
analytical continuation of $(s-1) \Phi (s,0)$ to the whole complex
plane. Moreover, for $\Real s$ from any bounded interval, the
estimate
$$
\Phi^{*} (s) = O(1)
$$
holds, where constant in $O$~- symbol is absolute. For $y>0$,
integration by parts twice gives
$$
\Phi(s,y) = O(|s|^2 y^{-2}).
$$

Square out the modulus of the sum in the expression (\ref{eq7}) for
$G(x)$, we get
\begin{equation*}
\begin{split}
& G (x)  =\sum\limits_{\nu_1,\nu_2, \nu_3, \nu_4} \frac{\beta(\nu_1)
\overline{\beta(\nu_2)\beta(\nu_3)}\beta(\nu_4)}{\nu_2\nu_4} \left(
\sum\limits_{n=1}^{+\infty} r(n) \overline{r\left( \frac{n
m_1}{m_2}\right)} \exp \left(
-\frac{4\pi}{\sqrt{D}}\frac{n m_1}{Q} x \sin\delta\right) \right.\\
& \left. + 2 \Real \sum\limits_{l\ge 1} \exp \left( -\frac{2\pi
i}{\sqrt{D}} \frac{l}{Q} x\cos\delta\right)
\sum\limits_{n=1}^{+\infty} r(n) \overline{r\left( \frac{n m_1 +l
}{m_2}\right)} \exp \left( -\frac{4\pi}{\sqrt{D}} \frac{(n m_1 +
l/2)}{Q} x \sin\delta \right) \right),
\end{split}
\end{equation*}
where the following notion were used
\begin{equation*}
\begin{split}
q = (\nu_1 \nu_4, \nu_2 \nu_3), \quad Q= \frac{\nu_2 \nu_4}{q},
\quad m_1 = \frac{\nu_1 \nu_4}{q} , \quad m_2 = \frac{\nu_2
\nu_3}{q}.
\end{split}
\end{equation*}
By Mellin's transform formula for Euler Gamma-function, $e^{-x} =
\frac{1}{2\pi i} \int\limits_{2-i\infty}^{2+i\infty}\Gamma(s) x^{-s}
ds$, we find
$$
J(x, \theta) \le J_1(x, \theta) + J_2(x, \theta),
$$
where
\begin{equation*}
\begin{split}
J_1(x, \theta) &= x^{1-\theta} \sum\limits_{\nu_1,\nu_2, \nu_3,
\nu_4}
\frac{\beta(\nu_1)\overline{\beta(\nu_2)\beta(\nu_3)}\beta(\nu_4)}{\nu_2\nu_4}
\times \\
& \times\left(  \frac{1}{2\pi i} \int\limits_{2-i\infty}^{2+i\infty}
D_{m_1, m_2} (s,0) \Gamma(s) B^{-s} \Phi (s+\theta, 0) ds \right),
\\
J_2(x, \theta) &= 2 x^{1-\theta} \sum\limits_{\nu_1,\nu_2, \nu_3,
\nu_4}
\frac{\beta(\nu_1)\overline{\beta(\nu_2)\beta(\nu_3)}\beta(\nu_4)}{\nu_2\nu_4}
\times \\
& \times \left(  \Real \sum\limits_{l\ge 1} \frac{1}{2\pi i}
\int\limits_{2-i\infty}^{2+i\infty} D_{m_1, m_2} (s,l) \Gamma(s)
B^{-s} \Phi (s+\theta, B_1 l)ds \right), \\
D_{m_1, m_2} &(s,l) = \sum\limits_{n} \frac{r(n) \overline{r\left(
\frac{m_1 n+l}{m_2}\right)}}{(m_1 n+l/2)^s},
\end{split}
\end{equation*}
\begin{equation*}
B = \dfrac{4\pi x\sin\delta}{\sqrt{D}Q}, \quad  B_1 = \dfrac{x
\cos\delta}{\sqrt{D}Q}.
\end{equation*}

For a given positive integer $m$, denote by $M_1 (m)$ the subset of
positive integers constituted by $1$ and by all numbers with the
same set of prime divisors as $m$ has. Since $(m_1, m_2)=1$, we have
\begin{equation*}
\begin{split}
D_{m_1, m_2} (s) &= D_{m_1, m_2} (s, 0) = \frac{1}{(m_1
m_2)^s}\sum\limits_{n=1}^{+\infty} \frac{r(n m_2) \overline{r(n m_1)}}{n^s} \\
&= \frac{1}{(m_1 m_2)^s} K(m_1, s) \overline{K}(m_2, s) D(s),
\end{split}
\end{equation*}
where
\begin{equation*}
\begin{split}
D(s) &=  \sum\limits_{n =1}^{+\infty} \frac{|r(n)|^{2}}{n^{s}}, \\
K(m,s) &=  \prod\limits_{p| m}
\left( 1+ \frac{|r(p)|^{2}}{p^{s}} + \frac{|r(p^2)|^{2}}{p^{2s}}+\ldots \right)^{-1} \times \\
&\times \prod_{p^{\alpha} || m} \left( \overline{r(p^{\alpha})} +
\frac{\overline{r(p^{\alpha+1})} r(p)}{p^{s}} +
\frac{\overline{r(p^{\alpha+2})} r(p^2)}{p^{2s}}+\ldots \right)\\
& = \sum\limits_{k \in M_1(m)} \frac{\overline{r(m k)} r(k)}{k^{s}}
\left( \sum\limits_{k \in M_1(m)} \frac{|r(k)|^{2}}{k^{s}}
\right)^{-1}.
\end{split}
\end{equation*}
If a prime $p$ divides $m_j$, then $p>256$; this entails that, for
$\Real s \ge 1/2$,
$$
\left| 1+ \frac{|r(p)|^{2}}{p^{s}} +
\frac{|r(p^2)|^{2}}{p^{2s}}+\ldots \right| > 0
$$
and, therefore, for a fixed $m_j$, that $K(m_j,s)$ is an analytic
function in the region $\Real s \ge 1/2$. Moreover, for either
$m=m_1$ or $m=m_2$, in this region the following estimate is valid
\begin{equation}\label{eq11}
\begin{split}
|K(m,s)| & \le \prod\limits_{p|m} \left( 1-\frac{2^2}{256^{1/2}}-
\frac{3^2}{256} - \frac{4^2}{256^{3/2}}-\ldots \right)^{-1}\times \\
&\times \tau(m)\prod_{p^{\alpha} || m} \left( 1 +
\frac{\frac{\alpha+2}{\alpha +1}\cdot 2}{256^{1/2}} +
\frac{\frac{\alpha+3}{\alpha +1}\cdot 3}{256}+
\frac{\frac{\alpha+4}{\alpha +1}\cdot 4}{256^{3/2}}+\ldots \right) \\
& \le \tau (m) \prod_{p | m} \left( 1-\frac{1}{4}- \frac{1}{4^2} -
\frac{1}{4^3}-\ldots \right)^{-1}\left( 1 + \frac{\frac{3}{2}\cdot
2}{16} + \frac{\frac{4}{2}\cdot 3}{16^2}+ \frac{\frac{5}{2}\cdot
4}{16^3}+\ldots \right)
\\
& = \tau (m) \prod_{p | m}\frac32\left( \sum\limits_{n \ge 2}
\frac{n(n-1)}{2\cdot 16^{n-2}}\right) \le \tau (m) \prod_{p | m} 3
\le \tau(m) \tau_{3} (m) \le \tau_{6} (m),
\end{split}
\end{equation}
where $\tau_t (n) = \sum\limits_{n_1 \cdot \ldots \cdot n_t =n} 1$.
Also, $D(s)$ can be meromorphically continued to the whole complex
plane with the pole of first order at $s=1$ (see~\cite{Rankin1939}).
This gives a meromorphic continuation of $D_{m_1, m_2}(s)$ into
half-plane $\Real s
>1/2$ .

In order to estimate $J_1(x, \theta)$, move the contour of
integration in
\begin{equation*}
\begin{split}
\int\limits_{2-i\infty}^{2+i\infty} D_{m_1, m_2} (s) \Gamma(s)
B^{-s} \Phi (s+\theta, 0) ds = \int\limits_{2-i\infty}^{2+i\infty}
\frac{K(m_1, s) \overline{K}(m_2, s) D(s) \Gamma(s) \Phi^{*}
(s+\theta) (B m_1 m_2)^{-s}}{s+\theta -1} ds
\end{split}
\end{equation*}
to the line $\Real s = 2/3$. We, therefore, pass simple poles at
$s=1$, $s = 1-\theta$. Denoting as $\mathfrak{D}$ the residie of
$D(s)$ at $s=1$, we get
\begin{equation*}
\begin{split}
& \frac{1}{2\pi i} \int\limits_{2-i\infty}^{2+i\infty} D_{m_1, m_2}
(s) \Gamma(s) B^{-s} \Phi (s+\theta, 0) ds
=\frac{\mathfrak{D} K(m_1, 1) \overline{K}(m_2, 1) \Phi^{*} (1+\theta)}{B m_1 m_2\theta}\\
&+ \frac{K(m_1, 1-\theta) \overline{K}(m_2, 1-\theta) D(1-\theta)
\Gamma(1-\theta) \Phi^{*} (1)}{(B m_1 m_2)^{1-\theta}} + O \left(
\frac{(m_1 m_2)^{1/3}}{(B m_1 m_2)^{2/3}}\right)\\
& = \frac{c(\theta) K(m_1, 1) \overline{K}(m_2, 1)}{\theta x \sin
\delta} \frac{q}{\nu_1 \nu_3} + \frac{c'(\theta) K(m_1, 1-\theta)
\overline{K}(m_2, 1-\theta) D(1-\theta)}{(x\sin\delta)^{1-\theta}}
\left(
\frac{q}{\nu_1 \nu_3}\right)^{1-\theta} \\
& + O\left( \frac{1}{(x\delta)^{2/3}} \frac{(\nu_2 \nu_4)^{1/3}}{
(\nu_1 \nu_3 )^{1/3}} \right),
\end{split}
\end{equation*}
where $c(\theta), c'(\theta)$ are constants that depend on $\theta$
and, for $0 < \theta \le 1/4$, are bounded. Applying estimate
$D(1-\theta) \ll \theta^{-1}$ and lemma~\ref{l3}, we arrive at the
following relation:
\begin{equation*}
\begin{split}
J_1 (x, \theta) &\ll \frac{|S(0)|}{\theta x^{\theta} \delta} +
\frac{|S(\theta)|}{\theta \delta^{1-\theta}} +
\frac{x^{1/3}}{x^{\theta}\delta^{2/3}}\sum\limits_{\nu_1, \ldots,
\nu_4 \le X}
\frac{1}{(\nu_1 \nu_3)^{1/3} (\nu_2 \nu_4)^{2/3}} \\
&\ll \frac{\delta^{-1}}{\theta x^{\theta} \ln X} \left( 1+ (X^2
e^{\frac{1}{h_1}} \delta)^{\theta} + X^2 (\ln X) e^{\frac{1}{3 h_1}}
\delta^{1/3}\right) \ll \frac{\delta^{-1}}{\theta x^{\theta} \ln X},
\end{split}
\end{equation*}
since we used $X^3 e^{\frac{1}{3 h_1}} \delta^{1/3} \ll 1 $.

Now let us estimate the ``non-diagonal'' term $J_2(x, \theta)$. It
follows from corollary~\ref{cor1} to lemma~\ref{l4}, that we can
move the path of integration in the integral $J_2 (x, \theta)$ to
the line $\Real s = 11/12$. Using the statement of lemma~\ref{l4},
we obtain the estimate
\begin{equation*}
\begin{split}
J_2 (x, \theta) &\ll x^{1-\theta} \sum\limits_{\nu_1, \ldots, \nu_4
\le X} \frac{1}{\nu_2 \nu_4} \sum\limits_{l\ge 1} \frac{1}{B^{11/12}
B_1^2 l^2} \left( (m_1 m_2)^{-2/11} + l^{1/10} m_1^{-10/11}
m_2^{-1}\right) \\
&\ll \frac{\delta^{-11/12}}{x^{\theta}} \sum\limits_{\nu_1, \ldots,
\nu_4 \le X} \frac{Q^{2+11/12}}{\nu_2 \nu_4} (m_1 m_2)^{-2/11} \ll
\frac{\delta^{-11/12}}{x^{\theta}} \sum\limits_{\nu_1, \ldots, \nu_4
\le
X} \frac{Q^{1+11/12}}{q (m_1 m_2)^{2/11}}  \\
&\ll \frac{\delta^{-11/12}}{x^{\theta}} \sum\limits_{\nu_1, \ldots,
\nu_4 \le X} \frac{Q^{1+11/12}}{(\nu_1 \nu_2 \nu_3 \nu_4)^{2/11}}
\ll \frac{\delta^{-11/12}}{x^{\theta}} (X^{2})^{1+\frac{11}{12}}
(X^{\frac{9}{11}})^4 \\
&\ll \frac{\delta^{-11/12}}{x^{\theta}} X^{469/66} \ll
\frac{\delta^{-1}}{x^{\theta} \ln X},
\end{split}
\end{equation*}
since $X \le \delta^{-1/86}$. The lemma is proved.
\end{pf0}

\begin{center}
{\bf \S 5. Proof of lemma \ref{l3}: estimation of Selberg sums}
\end{center}

\begin{pf0}
Application of M\"obious inversion formula to $f$, that is
\begin{equation*}
f(q) = \sum\limits_{d|q} \sum\limits_{m|d} \mu(m)
f\left(\frac{d}{m}\right),
\end{equation*}
gives the identity
\begin{equation*}
\begin{split}
q^{1-\theta} K\left( \frac{\nu_1 \nu_4}{q}, 1-\theta\right)
\overline{K}\left( \frac{\nu_2 \nu_3}{q}, 1-\theta\right)&=
\sum\limits_{d|q}  \sum\limits_{m|d} \mu(m) \left( \frac{d}{m} \right)^{1-\theta} \times \\
&\times K\left( \frac{\nu_1 \nu_4 m}{d}, 1-\theta\right)
\overline{K}\left( \frac{\nu_2 \nu_3 m}{d}, 1-\theta\right).
\end{split}
\end{equation*}
Inserting this formula to the expression for $S(\theta)$, we find
that
\begin{equation*}
\begin{split}
S(\theta) &=  \sum\limits_{\nu_1,\ldots, \nu_4 \le X}
\frac{\beta(\nu_1) \overline{\beta(\nu_2)  \beta(\nu_3)}
\beta(\nu_4)}{ \nu_2 \nu_4 \left( \nu_1 \nu_3 \right)^{1-\theta}}
\sum\limits_{d|q}  \sum\limits_{m|d} \mu(m) \left( \frac{d}{m} \right)^{1-\theta} \times \\
&\times K\left( \frac{\nu_1 \nu_4 m}{d}, 1-\theta\right)
\overline{K}\left( \frac{\nu_2 \nu_3 m}{d}, 1-\theta\right) =
\sum\limits_{d \le X^2} \sum\limits_{m|d} \mu(m) \left( \frac{d}{m}
\right)^{1-\theta} |g (d,m)|^2,
\end{split}
\end{equation*}
where
$$
g(d,m) = \sum\limits_{\substack{\nu_1 \nu_4 \equiv 0 (\!\!\!\!\mod d) \\
\nu_j \le X}} \frac{\beta(\nu_1)\beta(\nu_4)}{\nu_1^{1-\theta}
\nu_4} K\left( \frac{\nu_1 \nu_4 m}{d}, 1-\theta\right).
$$
Further, represent numbers $\nu_j$ in the form $\nu_j = \delta_j
\nu'_j$, where $(\nu'_j,d)=1$, and all prime divisors of $\delta_j$
are those (coinside) of  $d$, i.e., $\delta_j \in M_1 (d)$. Then
\begin{equation*}
\begin{split}
g(d,m) = \sum\limits_{\substack{\delta_1 \delta_4 \equiv 0
(\!\!\!\!\mod d)
\\ \delta_j \in M_1(d)}} \frac{1}{(\delta_1^{1-\theta} \delta_4)}
\sum\limits_{\substack{\delta_j \nu_j \le X, \\ (\nu_j, d)=1}}
\frac{\beta(\delta_1 \nu_1) \beta(\delta_4 \nu_4)}{\nu_1^{1-\theta}
\nu_4} K\left( \frac{\delta_1 \delta_4 m}{d} \nu_1 \nu_4,
1-\theta\right).
\end{split}
\end{equation*}
From the definition (\ref{eq5}) of $\beta(\nu)$ it follows, that
\begin{equation*}
\begin{split}
g(d,m) &= \frac{1}{\ln^2 X} \sum\limits_{\substack{\delta_1 \delta_4
\equiv 0 (\!\!\!\!\mod d) \\ \delta_j \in M_1(d)}}
\frac{\alpha(\delta_1) \alpha(\delta_4)}{\delta_1^{1-\theta}
\delta_4}
K\left( \frac{\delta_1 \delta_4 m}{d}, 1-\theta\right)  \times \\
&\times \sum\limits_{\substack{\nu_j \le X/\delta_j \\ (\nu_j,d) =1
}} \frac{\alpha(\nu_1) \alpha(\nu_4)}{\nu_1^{1-\theta} \nu_4}
K\left( \nu_1 \nu_4, 1-\theta\right) \ln\frac{X}{\delta_1 \nu_1}
\ln\frac{X}{\delta_4 \nu_4}.
\end{split}
\end{equation*}
We apply the following identity, which is valid for two
multiplicative functions $f_1$ and $f_2$, that are non-zero only for
square free numbers:
$$
\sum\limits_{\substack{l_j \le y_j, \\ (l_j, a)=1, \\ (l_1, l_2)=1}}
f_1 (l_1) f_2 (l_2) f(l_1, l_2) = \sum\limits_{(r,a)=1} \mu(r)
f_1(r) f_2(r)  \sum\limits_{\substack{\lambda_1 r \le y_1,\\
(\lambda_1, ra)=1}} \sum\limits_{\substack{\lambda_2 r \le y_2,\\
(\lambda_2, ra)=1}} f_1 (\lambda_1)  f_2 (\lambda_2) f(r\lambda_1,
r\lambda_2).
$$
Therefore, we get
\begin{equation*}
\begin{split}
g(d,m) &= \frac{1}{\ln^2 X} \sum\limits_{\substack{\delta_1 \delta_4
\equiv 0 (\!\!\!\!\mod d) \\ \delta_j \in M_1(d)}}
\frac{\alpha(\delta_1) \alpha(\delta_4)}{\delta_1^{1-\theta}
\delta_4}
K\left( \frac{\delta_1 \delta_4 m}{d}, 1-\theta\right)  \times \\
&\times \sum\limits_{(n,d)=1} \frac{\alpha^2(n)}{n^{2-\theta}}
K(n^2, 1-\theta)
\sum\limits_{(r, dn)=1} \frac{\mu(r) \alpha^2(r) K^2(r, 1-\theta)}{r^{2-\theta}} \times \\
&\times \left( \sum\limits_{\substack{\lambda_1 \delta_1 nr \le X,
\\ (\lambda_1, dnr)=1}}
\frac{\alpha(\lambda_1) K(\lambda_1, 1-\theta)}{\lambda_1^{1-\theta}}  \ln\frac{X}{\delta_1 nr \lambda_1} \right) \times \\
&\times \left(\sum\limits_{\substack{\lambda_4 \delta_4 nr \le X, \\
(\lambda_4, dnr)=1}} \frac{\alpha(\lambda_4) K(\lambda_4,
1-\theta)}{\lambda_4}  \ln\frac{X}{\delta_4 nr \lambda_4} \right).
\end{split}
\end{equation*}
Further we show that, for $0\le \theta, \gamma \le \frac14$, $N\ge
1$, $X_1 \ge 1$, the following estimate holds
\begin{equation}\label{eq17}
S_{\theta}(X_1,\gamma,N) = \sum\limits_{\substack{\lambda \le X_1, \\
(\lambda, N)=1}} \frac{\alpha(\lambda) K(\lambda,
1-\theta)}{\lambda^{1-\gamma}} \ln\frac{X_1}{ \lambda} \ll
X_1^{\gamma} \sqrt{\ln (X_1 +2)} \prod\limits_{p|N} \left(
1+\frac{1}{p} \right)^{2}.
\end{equation}
The equality
\begin{equation*}
\begin{split}
&\prod_{p>256} \left( 1+\frac{|r(p)|}{p^s} +
\frac{|r(p)|^2}{4p^{2s}}\right) = \prod_{p>256} \left(
1+\frac{|r(p)|}{2p^s}\right)^2 \\
&= \left( \sum\limits_{n=1}^{+\infty} \frac{|\alpha(n)|}{n^s}
\right)^2 = \sum\limits_{n=1}^{+\infty} \frac{1}{n^s} \left(
\sum\limits_{n_1 n_2 =n} |\alpha (n_1) \alpha(n_2)|\right)
\end{split}
\end{equation*}
entails that the function $b(n) = \sum\limits_{n_1 n_2 =n} |\alpha
(n_1) \alpha(n_2)|$ is multiplicative, and that also due to the
inequality $|r(p)| < 4$ (where $p$ is any prime number), for $n\mid
d$,
$$
b(nd) \le b(d).
$$
Using this inequality and the estimate $|K(n, 1-\theta)| \le
\tau_{6} (n)$ (see (\ref{eq11})), we find that
\begin{equation*}
\begin{split}
\sum\limits_{\substack{\delta_1 \delta_4 \equiv 0 (\!\!\!\!\mod d) \\
\delta_j \in M_1(d)}} & \frac{|\alpha(\delta_1)
\alpha(\delta_4)|}{\delta_1 \delta_4} \left| K\left( \frac{\delta_1
\delta_4 m}{d}, 1-\theta\right)\right| \le \frac{1}{d}
\sum\limits_{n \mid d} \frac{\mu^2 (n)}{n} | K(nm, 1-\theta)|
\sum\limits_{\delta_1 \delta_4 = n d} |\alpha(\delta_1)
\alpha(\delta_4)|
\\ &\le \tau_{6} (m) \frac{b(d)}{d} \prod_{p\mid d} \left( 1+
\frac{\tau_{6} (p)}{p} \right).
\end{split}
\end{equation*}
From this and (\ref{eq17}) we get
\begin{equation*}
\begin{split}
S(\theta) &\ll X^{2\theta} (\ln X)^{-2} \sum_{d \le X^2} \frac{b^2
(d)}{d^{1+\theta}}  \prod\limits_{p|d} \left(
1+\frac{1}{p}\right)^{8} \left( 1+ \frac{\tau_{6}(p)}{p} \right)^2
\sum\limits_{m|d} \frac{\mu^2(m)
\tau_{6}^2(m)}{m^{1-\theta}} \\
&\le X^{2\theta} (\ln X)^{-2} \sum_{d \le X^2} \frac{b^2
(d)}{d^{1+\theta}}  \prod\limits_{p|d} \left(
1+\frac{1}{p}\right)^{8} \left( 1+ \frac{\tau_{6}(p)}{p} \right)^2
\left( 1+ \frac{\tau_{6}^2(p)}{p^{1-\theta}} \right).
\end{split}
\end{equation*}
Now if $p$ is sufficiently large, and $0\le \theta \le \frac14$,
then
$$
\left( 1+\frac{1}{p}\right)^{8} \left( 1+ \frac{\tau_{6}(p)}{p}
\right)^2 \left( 1+ \frac{\tau_{6}^2(p)}{p^{1-\theta}} \right) \le
1+\frac{1}{\sqrt{p}};
$$
whence,
$$
\prod\limits_{p|d} \left( 1+\frac{1}{p}\right)^{8} \left( 1+
\frac{\tau_{6}(p)}{p} \right)^2 \left( 1+
\frac{\tau_{6}^2(p)}{p^{1-\theta}} \right) \ll \prod\limits_{p|d}
\left( 1+\frac{1}{\sqrt{p}} \right) \le \sum\limits_{m|d}
\frac{1}{\sqrt{m}}.
$$
From this inequality we find:
\begin{equation*}
\begin{split}
S(\theta) &\ll X^{2\theta} (\ln X)^{-2} \sum_{d \le
X^2} \frac{b^2 (d)}{d^{1+\theta}} \sum_{m|d} \frac{1}{\sqrt{m}}\\
&\le X^{2\theta} (\ln X)^{-2} \sum_{m \le X^2} \frac{1}{\sqrt{m}}
\sum_{\substack{d \le X^2 \\ d\equiv
0 (\!\!\!\!\mod m)}} \frac{b^2(d)}{d^{1+\theta}} \\
&\le X^{2\theta} (\ln X)^{-2}\sum_{m \le X^2} \frac{b^2 (m)}{m^{3/2
+ \theta}} \sum_{d \le X^2}
\frac{b^2(d)}{d^{1+\theta}} \\
&\ll X^{2\theta} (\ln X)^{-2} \sum_{d \le X^2} \frac{b^2(d)}{d} \ll
X^{2\theta} (\ln X)^{-1}.
\end{split}
\end{equation*}
The last estimate in the previous formula is provided by the
following one
\begin{equation*}
\begin{split}
\sum\limits_{n=1}^{+\infty} \frac{b^2 (n)}{n^s} = &\prod_{p>256}
\left( 1+\frac{|r(p)|^2}{p^s}+\frac{|r(p)|^4}{16 p^{2s}}\right) =
\prod_{p>256} \left( 1+\frac{|r(p)|^2}{p^s}\right) \left(
1+\frac{|r(p)|^4}{16 p^{2s} (1+\frac{|r(p)|^2}{p^s})}\right) \\
&= \sum\limits_{n=1}^{+\infty} \frac{1}{n^s} \left( \sum\limits_{n_1
n_2 =n} b_1 (n_1) b_2 (n_2)\right),
\end{split}
\end{equation*}
where
$$
\sum\limits_{n=1}^{+\infty} \frac{b_1 (n)}{n^s} = \prod_{p>256}
\left( 1+\frac{|r(p)|^2}{p^s}\right), \quad
\sum\limits_{n=1}^{+\infty} \frac{b_2 (n)}{n^s} = \prod_{p>256}
\left( 1+\frac{|r(p)|^4}{16 p^{2s} (1+\frac{|r(p)|^2}{p^s})}\right)
$$
and, thus,
\begin{equation*}
\begin{split}
\sum_{d \le X^2} \frac{b^2(d)}{d} \le \sum\limits_{n_1 n_2 \le X^2}
\frac{b_1 (n_1) b_2 (n_2)}{n_1 n_2} \le  \prod_{p} \left(
1+\frac{|r(p)|^4}{16 p^{2}
(1+\frac{|r(p)|^2}{p})}\right)\sum\limits_{n_1 \le X^2}
\frac{|r(n_1)|^2}{n_1} \ll \ln X,
\end{split}
\end{equation*}
in view of the equality (\cite{Rankin1939})
$$
\sum\limits_{n\le x} |r(n)|^2 = cx + O(x^{3/5}), \quad c>0
$$
(see~\cite{Hecke1937} for $O$~- estimate of the present sum, which
is also enough to obtain the desired result).

We then left to show the truth of estimate (\ref{eq17}) for the sum
$S_{\theta}(X_1,\gamma,N)$. Without loss of generality, we may
assume that $X_1\ge 10$, and that, for every $p|N$, the condition
$p>256$ is fulfilled. For $\Real s
>1$, consider the generated function
\begin{equation*}
\begin{split}
h_{\theta, N} (s) &= \sum\limits_{(n, N)=1} \frac{\alpha(n) K(n,
1-\theta)}{n^{s}} = \prod\limits_{p\nmid N} \left( 1 +
\frac{\alpha(p) K(p, 1-\theta)}{p^{s}}\right) \\
&= \prod_{\substack{p>256
\\ (p, N)=1}} \left( 1 - \frac{r(p) K(p,
1-\theta)}{2p^{s}}\right)  = \prod_{p} \left( 1 +
\frac{|r(p)|^2}{p^{s}}+
\frac{|r(p^2)|^2}{p^{2s}}+ \ldots\right)^{-1/2}\times \\
& \times \prod_{\substack{p>256 \\
(p, N)=1}} \left( 1 + \frac{|r(p)|^2}{p^{s}}+
\frac{|r(p^2)|^2}{p^{2s}}+ \ldots\right)^{1/2}\left( 1 - \frac{r(p)
K(p, 1-\theta)}{2p^{s}}\right) \times \\
&\times \prod\limits_{p\mid (256!)N} \left( 1 +
\frac{|r(p)|^2}{p^{s}}+ \frac{|r(p^2)|^2}{p^{2s}}+
\ldots\right)^{1/2}.
\end{split}
\end{equation*}
If $0 \le \theta \le 1/4$ then the product
\begin{equation*}
\begin{split}
N_1(s,\theta) &=\prod_{\substack{p>256 \\
(p, N)=1}} \left( 1 + \frac{|r(p)|^2}{p^{s}}+
\frac{|r(p^2)|^2}{p^{2s}}+ \ldots\right)^{1/2}\left( 1 - \frac{r(p)
K(p, 1-\theta)}{2p^{s}}\right) \\
&= \prod\limits_{\substack{p>256 \\ (p, N)=1}} \left( 1 +
\frac{|r(p)|^2}{2 p^{s}} + O\left( \frac{1}{p^{2\sigma}}\right)
\right) \left( 1 -
\frac{|r(p)|^2}{2 p^{s}} + O\left( \frac{1}{p^{\sigma+1-\theta}}\right)\right)\\
&= \prod\limits_{\substack{p>256 \\ (p, N)=1}} \left( 1 + O\left(
\frac{1}{p^{2\sigma}}\right) + O\left(
\frac{1}{p^{\sigma+1-\theta}}\right)\right), \quad \text{где} \quad
s=\sigma +it,
\end{split}
\end{equation*}
defines an analytic function in the half-plane $\Real s
> 1/2$.
For the generated function we have the identity
\begin{equation*}
\begin{split}
h_{\theta, N} (s) = (D(s))^{-1/2} N_1(s, \theta) G_{N} (s),
\end{split}
\end{equation*}
where
$$
D(s) = \sum\limits_{n=1}^{+\infty} \frac{|r(n)|^2}{n^s},
$$
$$
G_{N} (s) = \prod\limits_{p\mid (256!)N} \left( 1 +
\frac{|r(p)|^2}{p^{s}}+ \frac{|r(p^2)|^2}{p^{2s}}+
\ldots\right)^{1/2}.
$$

By means of Perron summation formula we find that
\begin{equation*}
\begin{split}
S_{\theta}(X_1,\gamma,N) &= \frac{1}{2\pi i}
\int\limits_{1-i\infty}^{1+i\infty} h_{\theta, N} (s+1-\gamma)
\frac{X_1^{s}}{s^2} ds \\
&= \frac{X_1^{\gamma}}{2\pi i}\int\limits_{1-i\infty}^{1+i\infty}
\frac{X_1^{s}}{(s+\gamma)^2}\frac{N_1 (s+1, \theta) G_{N} (s+1)
}{\left( D(s+1) \right)^{1/2}} ds.
\end{split}
\end{equation*}
For $\Real s \ge 1$, the following estimates are valid:
$$
N_1(s, \theta) \ll 1, \quad \text{если} \quad 0 \le \theta \le 1/4,
$$
$$
G_{N} (s) \ll \prod_{p|N} \left( 1+\frac{1}{p}\right)^2 = G_N.
$$
Move the path of integration from the line $\Real s =1$ to the
contour constructed by the semicircle $\{|s|=(\ln X_1)^{-1}, \quad
\Real s \ge 0 \}$ and the two rays $\{ s= it, |t| \ge (\ln X_1)^{-1}
\}$. For $D(s)$ in the region $\Real s \ge 1$ we shall use an
estimate $|D(s)|^{-1} \ll |s-1|$. The integral over the semicircle
(which we denote as $K_1$) can be estimated in the following way
\begin{equation*}
\begin{split}
K_1 \ll X_1^{\gamma} G_{N} \frac{X_1^{(\ln X_1)^{-1}} (\ln^{-1}
X_1)^{1+1/2}}{(\ln^{-1} X_1)^2} \ll X_1^{\gamma} (\ln X_1)^{1/2}
G_{N}.
\end{split}
\end{equation*}
For the integral over the rays (which we denote as $K_2$), the
relation
$$
K_2 \ll X_1^{\gamma} G_{N} \int\limits_{(\ln X_1)^{-1}}^{+\infty}
\frac{t^{1/2}}{t^2} dt \ll X_1^{\gamma} (\ln X_1)^{1/2} G_{N}.
$$
holds. Thus, the estimate (\ref{eq17}) is obtained and, therefore,
the lemma is proved.
\end{pf0}


\end{document}